\documentclass[a4paper,12pt,leqno]{article}
\usepackage{latexsym}
\usepackage[all]{xy}

\usepackage{amssymb} % \mathbb
\usepackage{amsmath} % \mathbb
\usepackage{theorem}

\def\Z{{\mathbb{Z}}}% \Z == \mathbb{Z}
\def\K{{\mathbb{K}}}% \K == \mathbb{K}
\def\R{{\mathbb{R}}}% \R == \mathbb{R}
\def\C{{\mathbb{C}}}% \C == \mathbb{C}
\def\A{{\mathcal{A}}}% \A == \mathcal{A}
\def\B{{\mathcal{B}}}% \B == \mathcal{B}
\DeclareMathOperator{\codim}{codim}

\DeclareMathOperator{\Der}{Der}

\DeclareMathOperator{\pd}{pd}
\DeclareMathOperator{\depth}{depth}

\DeclareMathOperator{\Poin}{Poin}

\numberwithin{equation}{section}

\newcommand{\owari}{\hfill$\square$}

\theoremstyle{break}
\newtheorem{theorem}{Theorem}[section]
\newtheorem{prop}[theorem]{Proposition}
\newtheorem{cor}[theorem]{Corollary}
\newtheorem{lemma}[theorem]{Lemma}

\newtheorem{rem}[theorem]{Remark}

\newtheorem{problem}[theorem]{Problem}

\title{Characteristic polynomials, $\eta$-complexes and freeness of tame arrangements}
\author{Takuro Abe
}
\date{\today} 

\begin{document}

\maketitle

\begin{abstract}
We compare each coefficient of the reduced characteristic polynomial 
of a simple arrangement and that of its Ziegler restriction. 
As a consequence we can show that the former is not less than 
the latter in the category of tame arrangements. This is a generalization of 
Yoshinaga's freeness criterion for $3$-arrangements and also the recent result 
by the author and Yoshinaga. As a corollary, we can prove that 
a free arrangement is a minimal chamber arrangement, and 
we can give a freeness criterion in terms of chambers in the category of tame arrangements.
%Let us consider a family of tame arrangements the Ziegler restrictions of 
%which are the same. 
%Then we can show that 
%%if the dimension is four, or 
%%tameness conditions are satisfied, then 
%there is a lower bound of 
%the cardinality of chambers of 
%arrangements belonging to this family. The lower bound depends on  
%their Ziegler restriction and is determined in an algebraic way. 
%Also, in the four dimensional cases, we can show that the minimal 
%chamber arrangement in the family above is the same as the freeness of 
%arrangements, which is a higher dimensional version of Yoshinaga's 
%criterion.
\end{abstract}

\section{Introduction}
Let $\A$ be a central $\ell$-arrangement over an arbitrary field $\K$. 
Fix $H_0 \in \A$ and $(\A'',m)$ the Ziegler restriction of $\A$ 
onto $H_0$. Let $d\A$ be the deconing of $\A$ with respect to $H_0$. For details of 
a notation in this section, see the next section.

Let us put a reduced characteristic polynomial of $\A$, which is combinatorial, as follows:
$$
\chi_0(\A,t)=\chi(d\A,t)=\sum_{i=0}^{\ell-1} (-1)^{\ell-1-i} b_{\ell-1-i} t^i.
$$
Also, let us put a characteristic polynomial of $(\A'',m)$, which is algebraic, as follows:
$$
\chi(\A'',m,t)=\sum_{i=0}^{\ell-1} (-1)^{\ell-1-i} \sigma_{\ell-1-i} t^i.
$$
It is known that $b_0=\sigma_0=1$ and $b_1=\sigma_1=|\A|-1=|m|$. Also, 
it is proved in \cite{AY} that $b_2 \ge \sigma_2$. Moreover, in \cite{AY}, 
the equality of $b_2$ and $\sigma_2$ is closely related to the freeness of $\A$. 
This is a generalization of Yoshinaga's freeness criterion for $3$-arrangements in \cite{Y}. 
After introducing a characteristic polynomial of multiarrangements in \cite{ATW}, 
Yoshinaga's criterion can be also understood 
in terms of the comparison of coefficients of characteristic polynomials, or 
minimality of chambers. 
%that to compare coefficients of 
%characteristic polynomial is important for the freeness criterion. 
Then a natural question is, what about $b_i$ and $\sigma_i$ for $i \ge 3$? 
The special case of this question is the relation between free arrangements and minimal chamber arrangements 
introduced in \cite{A}. To these problems, we can give an answer in 
the category of tame arrangements as follows:

\begin{theorem}
%Let $\A$ be a central $\ell$-arrangement. Fix 
%$H_0 \in \A$ and let $(\A'',m)$ be the Ziegler restriction 
%of $\A$ with respect to $H_0$. If $\A$ and $(\A'',m)$ are both 
%tame, then $(-1)^{\ell-1}\chi_0(\A,-1) \ge (-1)^{\ell-1}\chi(\A'',m,-1)  \ge 0$. In particular, it holds 
%that 
%$MC(\A'',m) \ge (-1)^{\ell-1} \chi(\A'',m,-1) \ge 0$ over the real number field. 
%%In particular, in the above condition, 
%%a free arrangement is a minimal chamber arrangement.
Let $\A$ be a central $\ell$-arrangement. Fix 
$H_0 \in \A$ and let $(\A'',m)$ be the Ziegler restriction 
of $\A$ with respect to $H_0$. 
%Put 
%\begin{eqnarray*}
%\chi_0(\A,t)&=&\sum_{i=0}^{\ell-1} (-1)^{\ell-1-i} b_{\ell-i-1} t^i\ (b_0=1),\\
%\chi(\A'',m,t)&=&\sum_{i=0}^{\ell-1} (-1)^{\ell-1-i} \sigma_{\ell-i-1} t^i\ (c_0=1).
%\end{eqnarray*}
If $\A$ and $(\A'',m)$ are both 
tame, then $b_i \ge \sigma_i \ge 0\ (i=0,1,\cdots,\ell-1)$ in the notation above.
\label{MCA}
\end{theorem}

Hence in the category of tame arrangements, we say that $\A$ is 
a \textbf{minimal chamber arrangement} (\textbf{MCA} for short) 
if $(-1)^{\ell-1}\chi_0(\A,-1)=(-1)^{\ell-1}\chi(\A'',m,-1)$. When  
$\ell=3$,  by Yoshinaga's criterion in \cite{Y}, we can define MCA,  
%Also, Yoshinaga's criterion in \cite{Y} 
and it holds that free arrangements and MCA are equivalent. 
As a corollary of Theorem \ref{MCA}, we can generalize this criterion and the relation between MCA
%minimal chamber 
%arrangement (MCA for short) 
and free arrangements as follows:

\begin{cor}
With the notation in Theorem \ref{MCA}, assume again that $\A$ and $(\A'',m)$ are tame. 
Then it holds that $(-1)^{\ell-1} \chi_0(\A,-1) \ge (-1)^{\ell-1} \chi(\A'',m,-1) \ge 0$. 
Moreover, $\A$ is free if and only if $(\A'',m)$ is free and $\chi_0(\A,-1) = \chi(\A'',m,-1)$. 
%In other words, in the tame category, MCA and free arrangements are equivalent.
\label{MCA2}
\end{cor}

Corollary \ref{MCA2} is also a generalization of Yoshinaga's criterion. In other words, 
if we fix a free multiarrangement $(\A'',m)$ and consider a family of arrangements $\A$ the Ziegler restriction of 
which are all $(\A'',m)$, then the freeness in this family is nothing but MCA in the tame 
category. Also, in the same category, a characteristic polynomial of the Ziegler 
restriction gives a lower bound of the value $(-1)^{\ell-1}\chi_0(\A,-1)$, or 
equivalently the cardinality of chambers over the real number field.  

The main tool for the proofs is the multi-version of the $\eta$-complex, originally 
introduced in \cite{ST}, developed 
in \cite{OT} and \cite{TY} for simple arrangements. In the proof, 
we also investigate several properties of this complex. 

The organization of this article is as follows. In section two
we introduce several definitions and results used in the rest of this  
article. In section three we 
develop several results for the proof.  
Mainly, we study several variants of the $\eta$-complexes. 
% \cite{ST}, \cite{OT} and \cite{TY}. 
In section four 
we prove Theorem \ref{MCA} and Corollary \ref{MCA2}. 
%The tools of the proof are 
%Poincar\'e polynomial, 
%the Solomon-Terao formula in \cite{ST} and 
%\cite{ATW}, and 
%In section four we prove Theorem \ref{nonnega} and Corollary \ref{Y4}.
%both 
%can be proved by an easy computation combined with the proposition 
%which appears in section two. 
%In section four we give rise to several 
%questions.
\medskip

\noindent
\textbf{Acknowledgements}. 
We thank Masahiko Yoshinaga for several comments to the draft of this 
article. The author is supported by JSPS Grants-in-Aid for Young Scientists
(B) No. 21740014.
%The author is supported by JSPS Grand-in-Aid for Young Scientists (B) No. 
%21740014. 
\medskip

\section{Preliminaries}
For the rest of this article everything is considered over an arbitrary field $\K$ and 
$V=\K^\ell$. For a general reference, see \cite{OT0}. 

Let $\A$ be an affine arrangement, i.e., 
a finite set of affine hyperplanes in $V$. 
An arrangement 
is called to be an $\ell$-arrangement if it is in $\K^\ell$. 
The \textbf{intersection lattice} $L(\A)$ is a set of subspaces of the form 
$\cap_{H \in \B} H$ with $\B \subset \A$. $L(\A)$ is a poset with the reverse inclusion order 
and the unique minimum element $V$. 
Define $L_i(\A)=\{X \in L(\A) \mid \codim_V X=i\}$. 
The \textbf{M\"obius function} $\mu:L(\A) \rightarrow 
\Z$ is defined by, $\mu(V)=1$, and by $\mu(X)=-\sum_{V \supset Y \supsetneq X} \mu(Y)\ (X \neq V)$. 
Then a \textbf{characteristic polynomial} $\chi(\A,t)$
% and 
%a \textbf{Poincar\'e polynomial} $\pi(\A,t)$ are 
is defined as follows:
%\begin{eqnarray*}
$$
\chi(\A,t)=\sum_{X \in L(\A)} \mu(X) t^{\dim X}.
%\\
%\pi(\A,t)&=&\sum_{X \in L(\A)} \mu(X) (-t)^{\codim X}.
$$
%\end{eqnarray*}
$\A$ is called to be \textbf{central} if $0 \in H\ (\forall H \in \A)$. 
Let $\alpha_H \in V^*$ be the defining form of $H \in \A$. 
If $\A$ is central, then $\chi(\A,t)$ 
%(resp: $\pi(\A,t))$ 
has $(t-1)$ 
%(resp: $(1+t)$) 
as a divisor. So define a \textbf{reduced
characteristic polynomial} $\chi_0(\A,t)$ 
% and \textbf{Poincar\'e polynomial} 
by 
%\begin{eqnarray*}
$$
\chi_0(\A,t):=\chi(\A,t)/(t-1).
%\\
%\pi_0(\A,t):&=&\pi(\A,t)/(1+t).
%\end{eqnarray*}
$$
Let $\A$ be a central $\ell$-arrangement.
$\A$ is called \textbf{essential} if $\cap_{H \in \A} H=\{0\}$. 
When $\A$ is a direct product of an essential arrangement 
$\B$ and an empty arrangement $\Phi$ (i.e., $\A \simeq \B \times \Phi)$, 
then $\B$ is called an 
\textbf{essentialization} of $\A$.
Now let us fix $H_ 0 \in \A$. 
Then the \textbf{deconing $d\A$} of $\A$ is defined as 
$\A \cap \{\alpha_{H_0}=1\}$, which is an $(\ell-1)$-affine arrangement. 
Note that $\chi_0(\A,t)=\chi(d\A,t)$. 
%and 
%$\pi_0(\A,t)=\pi(d\A,t)$. 
When the base field is $\R$, the set of connected components 
of $V \setminus \cup_{H \in \A} H$ is said to be 
\textbf{chambers}, and denoted by $C(\A)$.

\begin{rem}
It is well-known that $\pi(\A,t):=(-t)^\ell\chi(\A,-t^{-1})$ is equal to the topological 
Poincar\'e polynomial of $V \setminus \cup_{H\in \A} H$ when 
the base field is $\C$. Also, when the base field is $\R$, 
$(-1)^\ell \chi(\A,-1)$ is the number of chambers of the complement of 
hyperplanes, and 
$|\chi(\A,1)|$ the number of bounded chambers of that. 
Also, 
%it is well-known and easy to check that 
$|C(d\A)|=(-1)^{\ell-1} \chi_0(\A,-1)$.
%When the base field is arbitrary, there are no geometric or combinatorial 
%interpretation of these quantities, but we use the terminology that 
%$\pi(\A,1)$ is the number of chambers and so on.
\end{rem}

For the rest of this article we assume that $\A$ is a central 
$\ell$-arrangement. Let $S:=\mbox{Sym}^*(V^*)=\K[x_1,\ldots,x_\ell]$ be a coordinate ring of 
$V$. For the module of $S$-derivations $\Der S$, a \textbf{module of logarithmic vector fields of $\A$} is 
defined by  
$$
D(\A):=\{\theta \in \Der S \mid \alpha_H \mid \theta(\alpha_H)\ (\forall H \in \A)\}.
$$
In general $D(\A)$ is a reflexive module. When $D(\A)$ is a free $S$-module with 
homogeneous basis $\theta_1,\ldots, \theta_\ell$ of degrees $d_1,\ldots,d_\ell$, we say that 
$\A$ is \textbf{free} with \textbf{exponents} $\exp(\A)=(d_1,\ldots,d_\ell)$. 

A \textbf{multiplicity} is a map $m:\A \rightarrow \Z_{\ge 0}$ and a pair $(\A,m)$ 
is a \textbf{multiarrangement}. A \textbf{module of logarithmic vector fields of $(\A,m)$} is 
defined by 
$$
D(\A,m):=\{\theta \in \Der S \mid \alpha_H^{m(H)} \mid \theta(\alpha_H)\ (\forall H \in \A)\}.
$$
The freeness and exponents of a multiarrangement can be defined in the same manner. 
For $X \in L(\A)$, let $(\A_X,m_X)$ denote the \textbf{localization} of $(\A,m)$ defined by 
\begin{eqnarray*}
\A_X:&=&\{H \in \A \mid X \subset H\},\\
m_X:&=&m|_{\A_X}.
\end{eqnarray*}
%In particular, for a $2$-multiarrangement $(\A,m)$, we say that $(\A,m)$ is \textbf{balanced} 
%if $2m(H) \le |m|$ for any $H \in \A$. 

Multiarrangements appear naturally when we consider the restriction operation of a central 
arrangement. For a central arrangement $\A$ and $H_0 \in \A$, the \textbf{Ziegler restriction} 
$(\A'',m)$ with respect to $H_0$ is defined by 
\begin{eqnarray*}
\A'':&=&\{H \cap H_0 \mid H \in \A \setminus \{H_0\}\},\\
m(H \cap H_0):&=& |\{K \in \A \setminus \{H_0\} \mid K \cap H_0=H \cap H_0\}|.
\end{eqnarray*}

For 
the set of regular $p$-forms $\Omega^p_V$, a \textbf{module of logarithmic differential $p$-forms of 
$(\A,m)$} is defined as follows:
%we can define the dual object of $D(\A,m)$ as follows:
$$
\Omega^p(\A,m):=\{
\omega \in \displaystyle \frac{1}{Q(\A,m)} \Omega^p_V \mid 
(Q(\A,m)/\alpha_H^{m(H)})d \alpha_H \wedge \omega \in \Omega^{p+1}_V\ (\forall H \in \A)\},
$$
where 
$Q(\A,m):=\prod_{H \in \A} \alpha_H^{m(H)}$. See \cite{Z} for details of 
multiarrangements. By using this algebraic object, following 
\cite{ATW}, we can define 
a \textbf{characteristic polynomial of a multiarrangement} 
%and \textbf{Poincar\'e polynomial of a multiarrangement} 
as follows:
%\begin{eqnarray*}
$$
\chi(\A,m,t):=\lim_{x \to 1} \sum_{p=0}^{\ell} \mbox{Poin}(\Omega^p(\A,m),x)(t(1-x)-1)^p,
%\pi(\A,m,t):&=&(-t)^\ell\chi(\A,m,-t^{-1}).
$$
where $\mbox{Poin}(M,x):=\sum_{k \in \Z} \dim_\K M_k x^k$ is a Poincar\'e series of 
the $S$-graded module $M=\oplus_{k \in \Z}M_k$. 
%\end{eqnarray*}

\begin{rem}
Precisely, the definition of $\chi(\A,m,t)$ in the above is different from the original one 
in \cite{ATW}. In other words, the original definition was  
$$
\chi(\A,m,t):=(-1)^\ell\lim_{x \to 1} \sum_{p=0}^{\ell} \mbox{Poin}(D^p(\A,m),x)(t(x-1)-1)^p,
$$
where $D^p(\A,m)$ is a dual module of $\Omega^p(\A,m)$. 
The equality of these two definitions was proved in Remark 2.3 of \cite{AY}. So we use 
the definition by differential forms in this article.
%We will prove in Theorem \ref{domega} that 
%these definitions are equivalent. So we apply the definition by using 
%differential forms.
\end{rem}

Related to these characteristic polynomials, the following local-to-global formula 
is useful to compute each coefficient.

\begin{theorem}[\cite{ATW}, Theorem 3.3]
Put 
\begin{eqnarray*}
\chi(\A,m,t)&=&\sum_{i=0}^\ell (-1)^{\ell-i} \sigma_{\ell-i} t^{i},\\
\chi(\A_X,m_X,t)&=&t^{\ell-k}\sum_{i=0}^{k} (-1)^{k-i} \sigma_{k-i}^X t^{i}\ (X \in L_k(\A)).
\end{eqnarray*}
Then $\sigma_{k}=\sum_{X \in L_k(\A)} \sigma_{k}^X$.
\label{ATW}
\end{theorem} 

For a fixed $(\A'',m)$ where $\A''$ is a central $(\ell-1)$-arrangement, define 
$
F(\A'',m)$ to be the set of central $\ell$-arrangements the Ziegler restriction of which are all 
$(\A'',m)$. When it holds that 
$$
(-1)^{\ell-1} \chi_0(\A,-1) \ge (-1)^{\ell-1} \chi(\A'',m,-1) \ge 0
$$ 
for 
all $\A \in F(\A'',m)$, 
we say that $\A \in F(\A'',m)$ is a \textbf{minimal chamber arrangement (MCA} for short) 
if 
$$
(-1)^{\ell-1} \chi_0(\A,-1)=\min_{\B \in F(\A'',m)} (-1)^{\ell-1}\chi_0(\B,-1)=(-1)^{\ell-1}\chi(\A'',m,-1).
$$
We say that a multiarrangement $(\A,m)$ is \textbf{tame} if for 
a projective dimension $\pd_S \Omega^p(\A,m)$ of the $S$-module $\Omega^p(\A,m)$, it holds that 
$$
\mbox{pd}_S\Omega^p(\A,m) \le p\ (p=0,1,\ldots,\ell).
$$
%It is known that a tameness is a generic condition, 
For example, generic arrangements and free arrangements are tame, 
see \cite{RT} for 
details. Tame arrangements were introduced, first in \cite{OT} without names, and 
named in \cite{TY}. Recently, tame arrangements play important roles in several 
research areas of arrangements, see \cite{DSSWW}, \cite{DS} and \cite{Sc} for example. 

For $D(\A,m) \ni \theta$ and 
$\Omega^p(\A,m) \ni \omega=\sum g_{i_1\ldots i_p} dx_{i_1} \wedge \cdots \wedge 
dx_{i_p}$, define a \textbf{contraction} 
$$
\langle \theta,\omega \rangle:=
\sum (-1)^{j-1} \theta(x_{i_j})
g_{i_1\ldots i_p} dx_{i_1} \wedge \cdots \wedge dx_{i_{j-1}} \wedge dx_{i_{j+1}} \wedge \cdots \wedge
dx_{i_p}.
$$
If $\eta$ is a homogeneous $p$-form, then it holds that 
$$
\langle \theta, \eta \wedge \omega\rangle=\langle \theta,\eta \rangle \wedge \omega 
+(-1)^p \eta \wedge \langle \theta,\omega\rangle.
$$

The following is a generalized Yoshinaga's freeness criterion:

\begin{theorem}[\cite{AY}, Theorem 5.1]
Let $\A$ be an arrangement and $(\A'',m)$ the Ziegler restriction. 
Then $\A$ is free if and only if $(\A'',m)$ is free and 
$b_2=\sigma_2$ in the notation of the section one.
\label{AYfree}
\end{theorem}

%\begin{rem}
%Theorem \ref{Schulze} is generalized in \cite{AY}. In other words, 
%Theorem \ref{Schulze} holds true without the assumption that $\A$ is tame in 
%any dimensional vector space. Also, the freeness of $(\A'',m)$ and 
%$b_2=\sigma_2$ is sufficient. However, Theorem \ref{Schulze} is sufficient in this 
%article.
%\end{rem}

The following map, which is 
%a part of the result 
introduced in 
\cite{AY}, is important to prove Theorem \ref{MCA}. 

\begin{prop}[\cite{AY}]
Let $\A$ be a central $\ell$-arrangement, 
$H_0 \in \A$ and $(\A'',m)$ the Ziegler restriction. 
%Let $L(d\A)$ be the intersection lattice of the deconing of 
%$\A$ and $L(\A'')$ that of $\A''$. 
Then 
there is a well-defined map $\rho:L(d\A) \rightarrow L(\A'')$ which 
keeps inclusion orders and codimensions of each flat. Also, $\rho$ is 
compatible with localization operations. 
%Put 
%\begin{eqnarray*}
%\pi_0(\A,t)&=&1+b_1t+b_2t^2+\cdots+b_{\ell-1}t^{\ell-1},\\
%\pi(\A'',m,t)&=&1+\sigma_1t+\sigma_2t^2+\cdots+\sigma_{\ell-1}t^{\ell-1}.
%\end{eqnarray*}
%Then $b_2 \ge \sigma_2$.
\label{betti2}
\end{prop}

%\begin{rem}
%Theorem \ref{betti2} is the other direction of generalizing Yoshinaga's criterion. 
%Theorem \ref{betti2} generalizes the inequality between the second Betti 
%number and the coefficient of $t^2$ of the multi-characteristic polynomial. See 
%\cite{AY} for more details. 
%On the other hand, Theorem \ref{MCA} generalizes the lower bound 
%formula of betti numbers and chambers, or equivalently, $(-1)^{\ell-1}\chi_0(\A,t)$ in 
%terms of multi-characteristic polynomial. 
%\end{rem}

\section{Several complexes and their properties}

Put $\alpha:=\alpha_{H_0}=x_\ell \in S=\K[x_1,\ldots,x_\ell]$ and 
$S'=S/\alpha S=\K[x_1,\ldots,x_{\ell-1}]$ the coordinate ring of $H_0$. 
To prove Theorem \ref{MCA} we need some lemmas and propositions, mainly on 
$\eta$-complexes. 

\begin{rem}
In this section we do not use the tameness assumption.
\end{rem}
%Let $\K$ be an arbitrary field. 

\begin{lemma}
The $S$-morphism 
$$
\Omega^p(\A) \rightarrow \Omega^p(\A) \wedge \displaystyle \frac{d\alpha}{\alpha}
\rightarrow 0.
$$
is a splitting surjection. In particular, 
$\pd_S \Omega^p(\A) \ge \pd_S 
(\Omega^p(\A) \wedge \displaystyle \frac{d\alpha}{\alpha})$.
\label{lemma1}
\end{lemma}

\noindent
\textbf{Proof}. It suffices to show that the morphism has a section. Recall 
Proposition 4.86 in \cite{OT0}. Then the section is given by 
$$
\omega \wedge \displaystyle \frac{d \alpha}{\alpha} \mapsto 
(-1)^p\langle \theta_E, \omega \wedge \displaystyle \frac{d \alpha}{\alpha}\rangle,
$$
where $\langle, \rangle$ is a contraction. 
The inequality of projective dimensions follows from the long exact sequence of 
Ext's.  \owari
\medskip

\begin{rem}
Since the complex $(\Omega^*(\A),\wedge \displaystyle \frac{d \alpha}{\alpha})$ 
is 
exact (see \cite{OT0} for example), Lemma \ref{lemma1} shows that 
$$
\Omega^p(\A) \simeq (\Omega^{p-1}(\A) \wedge \displaystyle \frac{d\alpha}{\alpha}) \oplus 
(\Omega^p(\A) \wedge \displaystyle \frac{d\alpha}{\alpha}).
$$
\label{ds}
\end{rem}

Let $\mbox{res}:\Omega^p(\A) \rightarrow \Omega^p(\A'',m)$ be the residue map defined by 
$$
\sigma \wedge \displaystyle \frac{d \alpha}{\alpha}  + \delta \mapsto \delta|_{H_0},
$$
where $\sigma$ and $\delta$ are generated by $dx_1,\ldots,d x_{\ell-1}$. Note that 
the residue map factors through $\Omega^p(\A) \wedge \displaystyle \frac{d\alpha}{\alpha} 
\rightarrow \Omega(\A'',m)$. Let $M^p \subset \Omega^p(\A'',m)$ denote the image 
of the residue map and $C^p$ its cokernel:
$$
0 \rightarrow M^p \rightarrow \Omega^p(\A'',m) \rightarrow C^p \rightarrow 0.
$$ 

\begin{lemma}
The sequence 
$$
0 \rightarrow \Omega^p(\A) \wedge \displaystyle \frac{d\alpha}{\alpha} \rightarrow 
\Omega^p(\A) \wedge \displaystyle \frac{d\alpha}{\alpha}
\rightarrow M^p \rightarrow 0
$$
is exact, where the second arrow is the product of $\alpha$ and 
the third arrow is the residue map. In particular, 
$\pd_{S'} M^p \le \pd_S \Omega^p(\A)$. 
\label{lemma2}
\end{lemma}

\noindent
\textbf{Proof}. 
$$
\mbox{res}(\delta \wedge \displaystyle \frac{d \alpha}{\alpha}  )=\delta|_{H_0}=0 \iff
\alpha \mid \delta.
$$
Hence the exactness follows immediately. Let us prove the inequality. Since the action of 
$S$ to $M^p$ factors through $S'=S/\alpha S$, it follows that $\depth_S M^p=\depth_{S'} M^p$. 
Hence Auslander-Buchsbaum formula shows that 
$\pd_{S'} M^p+1=\pd_S M^p$. Also, the long exact sequence shows that 
$\pd_S (\Omega^p(\A) \wedge \displaystyle \frac{d \alpha}{\alpha}) +1 \ge \pd_S M^p$. Combining 
this with Lemma \ref{lemma1} gives $\pd_{S'} M^p \le \pd_S \Omega^p(\A)$. \owari
\medskip

Next let us consider the $\eta$-complex, see \cite{OT0} for details. 
It is the complex $(\Omega^*(\A),\wedge \eta)$, where $\eta$ is some generic regular $1$-form and 
the boundary map is given by $\wedge \eta$. This is of course a complex, and we can 
define the cohomology group $H^p(\Omega^*(\A))$. Let $\overline{\eta}:=
\eta|_{H_0}$. Since the wedge product is commutative with the 
inclusion $M^p \rightarrow \Omega^p(\A'',m)$ and $\eta$ is regular, 
the wedge product of $\overline{\eta}$ is closed in $\Omega^p(\A'',m)$. In other words, 
the boundary map $\wedge \overline{\eta}:\Omega^p(\A'',m) \rightarrow \Omega^{p+1}(\A'',m)$ is induced for 
$p=0,1,\ldots,\ell-1$.  
Since $M^p$ 
and $C^p$ are the surjective images from these two differential modules, 
we can define not only the complexes 
$(M^*, \wedge \overline{\eta}),\ (\Omega^*(\A'',m),\wedge \overline{\eta})$ and 
$(C^*,\wedge \overline{\eta})$ but also the cohomology groups 
$H^p(M^*), H^p(\Omega^*(\A'',m))$ and $H^p(C^*)$.
%, where 
%$C^p:=\coker (M^p \rightarrow \Omega^p(\A'',m))$. 

\begin{prop}
For an integer $d \ge 0$, 
there exists a regular generic $1$-form $\eta$ of homogeneous degree $d$ such that 
all cohomology groups of both complexes $(\Omega^*(\A),\wedge \eta)$ and 
$(\Omega^*(\A'',m), \wedge \overline{\eta})$ are finite dimensional.
\label{finite}
\end{prop}

\noindent
\textbf{Proof}. The proof is similar to that in \cite{OT0} with a slight modification for 
multiarrangements. 
Let $S^X$ be the coordinate ring of $X \in L(\A)$. 
Let $r_{Y,X}:S^Y \rightarrow S^X$ be the quotient map for flats $X \subset Y$ in $L(\A)$, and 
$\Omega^1[X]^0_d$ the set of regular $1$-forms of degree $d$ over $X$ which 
vanish only at the origin. It is well-known that such forms are generic in each 
vector spaces. 
Also, we can canonically extend $r_{Y,X}$ to that from the set of 
differential forms over $Y$ to those over $X$. 
Now put $N^X_d:=r_{V,X}^{-1} (\Omega^1[X]^0_d)$ and 
$$
N_d:=\bigcap_{X \in L(\A),\ \dim X>0} N^X_d.
$$
Since $r_{V,X}$ is continuous, $N_d$ is a non-empty open set. Take an arbitrary $\eta \in N_d$. 
Then Proposition 4.91 in \cite{OT0} shows $\dim_\K H^p(\Omega^*(\A))< \infty$. 
%
%Let us prove that, for the ideal 
%$$
%I(\eta):=\{\langle \theta, \eta \rangle \mid \theta \in D(\A)\},
%$$
%the radical of $I(\eta)$ contains the irrelevant ideal of $S$.
%% which 
%%is proved in \cite{OT}. 
%By a change of base fields we may assume that $\K$ is an algebraically closed field. 
%Then by Hilbert's Nullstellensatz, it suffices to show that 
%the zero locus $Z(I(\eta))$ of the ideal $I(\eta)$ is contained in the origin. 
%Take $v \in V \setminus \{0\}$ and put 
%$X:=\cap_{v\in H \in \A} H$. Assume that $v \in Z(I(\eta))$. Choose a basis 
%$x_1,\ldots,x_\ell$ for $V^*$ in such a way that 
%$X=\{x_{k+1}=\cdots=x_\ell=0\}$. Put 
%$\A_1:=\{H \in \A \mid X \not \subset H\}$ and 
%$Q_1:=Q(\A_1)$. By definition 
%$Q_1 \partial_{x_i} \in D(\A)$ for $i=1,\ldots,k$. Write 
%$\eta=f_1dx_1+\cdots+f_\ell dx_\ell$ with $f_i \in S$. 
%Then $\langle Q_1 \partial_{x_i},\eta\rangle=Q_1f_i \in I(\eta)$. By definition 
%$Q_1(v) \neq 0$. Hence $f_1(v)=\cdots=f_k(v)=0$. Recalling that 
%$r_{V,X}(\eta)=\overline{f}_1dx_1+\cdots+\overline{f}_k dx_k$ vanishes only at the origin, 
%there exists some $i,\ 1 \le i \le k$ such that ${f}_i(v) \neq 0$, which is a contradiction. 
%
Next, let us prove the multi-case. First, let us prove that, 
for the ideal 
$$
I(\overline{\eta}):=\{\langle \overline{\theta}, \overline{\eta}\rangle \mid \overline{\theta} \in D(\A'',m)\} 
\subset S',
$$
the radical of $I(\overline{\eta})$ contains the irrelevant ideal of $S'$. 
It suffices to show that 
the zero locus $Z(I(\overline{\eta}))$ of the ideal $I(\overline{\eta})$ is contained in the origin. 
Take $v \in H_0 \setminus \{0\}$ and put 
$X_0:=\cap_{v\in H' \in \A''} H'$. Assume that $v \in Z(I(\overline{\eta}))$. Choose a basis 
$x_1,\ldots,x_{\ell-1}$ for $H_0^*$ in such a way that 
$X_0=\{x_{k+1}=\cdots=x_{\ell-1}=0\}$. Put 
$\A_1'':=\{H' \in \A'' \mid X_0 \not \subset H'\}$ and 
$Q_1':=Q(\A_1'',m|_{\A_1''})$. By definition 
$Q_1' \partial_{x_i} \in D(\A'',m)$ for $i=1,\ldots,k$. Write 
$\overline{\eta}=f_1dx_1+\cdots+f_{\ell-1} dx_{\ell-1}$ with $f_i \in S'$. 
Then $\langle Q_1' \partial_{x_i},\overline{\eta}\rangle=Q_1'f_i \in I(\overline{\eta})$. By definition 
$Q_1'(v) \neq 0$. Hence $f_1(v)=\cdots=f_k(v)=0$. Recalling that 
$r_{H_0,X_0}(\overline{\eta})=\overline{f}_1dx_1+\cdots+\overline{f}_k dx_k$ vanishes only at the origin, 
there exists some $i,\ 1 \le i \le k$ such that ${f}_i(v) \neq 0$, which is a contradiction.

%Now we can show that there exists a common regular generic $1$-form $\eta$ of degree $d$ such that 
%$I(\eta)$ and $I(\overline{\eta})$ annihilate each of corresponding cohomology groups respectively.
%Let $\eta \in N_d$ and $H^p$ denote the $p$-th cohomology of the $\eta$-complex 
%$(\Omega^*(\A),\eta)$. 
%Let $\omega \in \Omega^p(\A)$ be a cocycle of this complex and take $\theta \in D(\A)$. Then 
%$$
%0=\langle \theta, \eta \wedge \omega \rangle
%=\langle \theta, \eta \rangle \omega -
%\eta \wedge \langle \theta, \omega \rangle.
%$$
%Hence $I(\eta)$ annihilates $H^p$. 

Second, let $\eta \in N_d$ and $\overline{H}^p$ denote the $p$-th cohomology of the $\overline{\eta}$-complex 
$(\Omega^*(\A'',m),\wedge \overline{\eta})$ where $\overline{\eta}=r_{V,H_0} (\eta)$. First, note that 
$\overline{\eta} \neq 0$ since $\overline{\eta}$ is chosen in such a way that 
it only vanishes at the origin. 
Before the proof, let us show the following easy but important lemma.

\begin{lemma}
For $\theta \in D(\A'',m)$ and $\omega \in \Omega^p(\A'',m)$, it holds that 
$\langle \theta,\omega\rangle \in \Omega^{p-1}(\A'',m)$.
\label{contract}
\end{lemma}

\noindent
\textbf{Proof}. 
For 
$H \in \A''$, it holds that 
$$
\langle \theta, d \alpha_H \wedge \omega\rangle=
\langle \theta, d \alpha_H \rangle \omega -d\alpha_H \wedge \langle \theta,\omega \rangle
=\theta(\alpha_H) \omega -d\alpha_H \wedge \langle \theta,\omega \rangle.
$$
Since $\alpha^{m(H)} \mid \theta(\alpha_H)$ and $d\alpha_H \wedge \omega$ is regular along $H$,  
$d \alpha_H \wedge \langle \theta,\omega \rangle$ is regular along $H$. 
Since $Q(\A'',m) \langle \theta,\omega\rangle$ is regular, 
$\langle \theta,\omega\rangle\in \Omega^{p-1}(\A'',m)$. \owari
\medskip

\noindent
\textbf{Proof of Proposition \ref{finite}, continued}. 
Now let $\overline{\omega} \in \Omega^p(\A'',m)$ be a cocycle of this complex and take $\overline{\theta} \in D(\A,m)$.  
Then 
$$
0=\langle \overline{\theta}, \overline{\eta} \wedge \overline{\omega} \rangle
=\langle \overline{\theta}, \overline{\eta} \rangle \overline{\omega} -
\overline{\eta} \wedge \langle \overline{\theta}, \overline{\omega} \rangle.
$$
Hence $I(\overline{\eta})$ annihilates $\overline{H}^p$, 
%see the arguments in \cite{OT}, 
which makes the cohomology group finite dimensional. 
\owari
\medskip

\begin{rem}
Proposition \ref{finite} shows that we can choose the regular $1$-form 
$\eta$ in the proposition such that the dimensions of all cohomologies of the $\eta$-complex are 
finite for all $m:\A \rightarrow \Z_{\ge 0}$. In other words, 
such a $1$-form $\eta$ depends only on $\A$, independent of $m$. 
\end{rem}

%Then it is easy to show that $\wedge \eta$ and $\wedge \overline{\eta}$ induces 
%the complex structure not only on $\Omega^*(\A)$ and $\Omega^*(\A'',m)$ but also 
%$M^*$ and $C^*$. 

\begin{cor}
In the same notation, $H^p(M^*)$ is also finite dimensional.
\label{finite1}
\end{cor}

\noindent
\textbf{Proof}. 
First, by the exact sequence in Lemma \ref{lemma1}, Proposition \ref{finite} and 
Remark \ref{ds}, 
it holds that $H^p(\Omega^*(\A) \wedge d \alpha/\alpha)$ is of finite 
dimensional. So the exact sequence in Lemma \ref{lemma2} shows that $H^p(M^*)$ is 
all finite 
dimensional. \owari
\medskip

\begin{cor}
In the same notation, $H^p(C^*)$ is also finite dimensional.
\label{finite2}
\end{cor}

\noindent
\textbf{Proof}. Apply 
Proposition \ref{finite} and Corollary \ref{finite1}
to the cohomology long exact sequence of 
$$
0 \rightarrow M^p \rightarrow \Omega^p(\A'',m) \rightarrow C^p \rightarrow 0
$$
which commutes with $\wedge \overline{\eta}$.  \owari
\medskip

Before the next proposition, let us recall the fact that 
$C^0=C^{\ell-1}=0$. We follow the proof in \cite{Sc}. Since 
$\Omega^0(\A)=S$ and 
$\Omega^0(\A'',m)=S'$, it follows that $C^0=0$. 
Also, since the complex $(\Omega^*(\A),\wedge \displaystyle \frac{d\alpha}{\alpha})$ is 
exact, it follows that 
$$
\Omega^{\ell-1}(\A) \wedge \displaystyle \frac{d\alpha}{\alpha}
=\Omega^\ell(\A)=S/Q(\A)dx_1 \wedge \cdots \wedge dx_\ell.
$$
So 
$$
\Omega^{\ell-1}(\A'',m)=S'/Q(\A'',m) dx_1\wedge \cdots \wedge dx_{\ell-1}
$$
implies that $C^{\ell-1}=0$.

\begin{prop}
$\sum_{p=0}^{\ell-1} \mbox{Poin}(C^p,x)(t(1-x)-1)^p \in 
\R[x,x^{-1},t]$ and 
$\sum_{p=0}^{\ell-1} \mbox{Poin}(M^p,x)(t(1-x)-1)^p \in 
\R[x,x^{-1},t]$, i.e., there are no poles along $x=1$. 
\label{poly}
\end{prop}

\noindent
\textbf{Proof}. 
Apply the same proof as Proposition 4.133 in \cite{OT0} 
combined with Propositions \ref{finite}, Corollaries \ref{finite1} and \ref{finite2}. \owari
%First let us generalize Proposition \ref{finite} for 
%regular generic $d$-forms, which can be proved by the same arguments 
%in \cite{ST} and \cite{ATW}. Then also the same arguments 
%show that the finiteness of cohomology groups of these 
%degree $d$-complexes of $C^*$. 

%省いたパート：ここから
%Since the proofs are the same, we only show that 
%$\sum_{p=0}^{\ell-1} \mbox{Poin}(C^p,x)(t(1-x)-1)^p \in 
%\R[x,x^{-1},t]$. 
%First, take minimum integers $n$ and $m$ such that the series 
%$$
%P(x,t):=x^n (1-x)^m\sum_{p=0}^{\ell-1} \mbox{Poin}(C^p,x)(t(1-x)-1)^p
%$$
%is a polynomial both in $x$ and $t$. Note that $m \le \ell-1$. Hence for $d \ge 0$, 
%$$
%P(x,(1-x^{-d})/(1-x))=x^n (1-x)^m \sum_{p=0}^{\ell-1} \mbox{Poin}(C^p,x)(-x^{-d})^p.
%$$
%Then Proposition \ref{finite} allows us to choose a regular generic 
%$1$-form $\eta_d$ of degree $d$ which makes cohomology groups of the $\eta_d$-complex finite 
%dimensional. Then 
%$$
%\sum_{p=0}^{\ell-1} \mbox{Poin}(C^p,x)(-x^{-d})^p
%=\sum_{p=0}^{\ell-1} \mbox{Poin}H^p(C^*)(-x^{-d})^p
%$$
%is just an alternating sum of each $\eta_d$-complexes and 
%has a finite value when $x=1$. Assume that $m>0$. Then 
%$$
%P(x,(1-x^{-d})/(1-x))=x^n (1-x)^m \sum_{p=0}^{\ell-1} \mbox{Poin}(C^p,x)(-x^{-d})^p
%$$
%is equal to zero when $x=1$. Hence for infinitely many $d \in \Z_{>0}$, it holds that 
%$P(1,-d)=0$. Since $P(1,t)$ is a polynomial, it follows that $P(1,t)=0$. This implies that 
%$\sum_{p=0}^{\ell-1} \mbox{Poin}(C^p,x)(t(1-x)-1)^p$ has no poles along $x=1$. 
%Equivalently, $m=0$, which is a contradiction. \owari
%省いたパート：ここまで
\medskip
%Then taking the alternating sum of the 
%Poincar\'e polynomials gives us that 
%\begin{eqnarray*}
%\sum_{p=0}^{\ell-1} \mbox{Poin}(C^p,x)(-x^d)^p=
%\sum_{p=0}^{\ell-1} \mbox{Poin}(H^p(C^*),x)(-x^d)^p
%\end{eqnarray*}
%is a polynomial in $x$ for infinitely many $d \in \Z$. 
%
%Now recall that $\sum_{p=0}^{\ell-1} (x-1)^{\ell-1} \mbox{Poin}(C^*,x)y^p$ is a 
%polynomial. Hence we can expand this into 
%$$
%\sum_{i,j \ge 0} c_{ij}(x-1)^i(y+1)^j\ (c_{ij} \in \Q).
%$$
%Then the same arguments as in Proposition 5.2 in \cite{ST} completes the proof. \owari
%\medskip

The following is useful to prove Theorem \ref{MCA}.

\begin{theorem}[Theorem 5.8, \cite{OT}]
Let  $S=\K[x_1,\ldots,x_\ell]$  and $F^*=(0 \rightarrow F^0 \rightarrow 
F^1 \rightarrow \cdots \rightarrow F^\ell \rightarrow 0)$ be a complex of finite $S$-modules such that 
every morphism is $S$-linear and that 
every cohomology group is finite dimensional. 
If a nonnegative integer $q$ satisfies 
$$
\pd_S F^p < \ell +p-q
$$
for all $p$, then $H^q(F^*)=0.$
\label{TY72}
\end{theorem}

\section{Proof of Theorem \ref{MCA} and Corollary \ref{MCA2}}

In this section we prove the main results of this article. Recall that 
we have not yet used the tameness assumption in this article. In this section we apply 
it. 
\medskip

\noindent
\textbf{Proof of Theorem \ref{MCA}}. 
Let us prove Theorem \ref{MCA} by induction on the dimension $\ell$ and $i$ in the setup of 
Theorem \ref{MCA}. For $i=0$, $b_0=\sigma_0=1$. For $i=1$, 
%the fact that 
$b_1=\sigma_1=|\A|-1=|m|$.  
%implies the statement. 
So Theorem \ref{MCA} holds when $\ell \le 1$. For $\ell=2$, as we see in 
the section one, the statement is nothing but Yoshinaga's criterion (see also 
\cite{AY}). Assume that $\ell \ge 3$ and 
$i<\ell-1$. Recall the map $\rho:L(d\A) \rightarrow L(\A'')$ in Proposition \ref{betti2} and 
put 
$$
b_i =\sum_{X \in L_i(\A'')}b_i^X,
$$
where $b_i^X$ is the sum of absolute values of $\mu(Y)$ with 
$Y \in L_i(d\A)$ and $\rho(Y)=X$. Also, let 
$
\sigma_i^X
$ be the absolute value of the constant term 
of $\chi_{red}(\A''_X,m_X,t):=\chi(\A''_X,m_X,t)/t^{\ell-1-i}$ for 
$X \in L_i(\A'')$. 
%Also put $\chi_{0,red}(\A_X,t):=\chi_0(\A_X,t)/t^{\ell-1-\codim X}$. 
%, where 
%${(\A'')}_X^e$ is the essentialization of the localization of $\A''$ at $X$. 
Then
%local-to-global formula in \cite{ATW} and 
Theorem \ref{ATW} and Proposition \ref{betti2} imply that 
$$
%b_i=\sum_{X \in L_i(\A'')}b_i^X
%=\sum_{X \in L_i(\A'')}|\chi_{0,red}(\A_X,0)|
%$$
%and 
%$$
\sigma_i=\sum_{X \in L_i(\A'')}\sigma_i^X=\sum_{X \in L_i(\A'')}|\chi_{red}(\A''_X,m_X,0)|.
$$
%Note that 
%$$
%b_i=
%$$
Now recall the tameness condition on $\A$ and $(\A'',m)$, which we have not yet used. 
By definition of the tameness and the fact that the localization is exact, 
it holds that $\A_X$ and $(\A''_X,m_X)$ are also tame. However, $\A_X$ and 
$(\A''_X,m_X)$ are both $\ell$-multiarrangements. Hence to apply the induction hypothesis, 
we need the following lemma (see also \cite{DS}):

\begin{lemma}
Let $\A=\A_1 \times \A_2$ be an $\ell$-arrangement which decomposes into 
the product of a $d$-arrangement $\A_1$ in $V_1$ and an $(\ell-d)$-arrangement $\A_2$ in $V_2$. 
For $m:\A \rightarrow \Z_{\ge 0}$, let $m_i\ (i=1,2)$ denote 
the restriction of $m$ onto $\A_i$. Let $S_i\ (i=1,2)$ denote 
the coordinate ring of $V_i$. Hence $V_1 \otimes_\K V_2 = V$ and $S = S_1 \otimes_\K S_2$. 
Then 
$\pd_S \Omega^p(\A,m) \ge \pd_{S_i} \Omega^p(\A_i,m_i)$.
\label{pd}
\end{lemma}

\noindent
\textbf{Proof}. 
It suffices to show when $i=1$. 
Since 
$$
\Omega^p(\A,m)=\oplus_{q+r=p} \Omega^q(\A_1,m_1) \otimes_\K \Omega^r(\A_2,m_2),
$$
it holds that $\pd_S \Omega^p(\A,m) \ge \pd_S S \cdot \Omega^p(\A_1,m_1)$. 
Let us prove $\pd_S S \cdot \Omega^p(\A_1,m_1) \ge \pd_{S_1} \Omega^p(\A_1,m_1)$. 

Note that $S\cdot \Omega^p(\A_1,m_1) \simeq \Omega^p(\A_1,m_1) \otimes_{\K} S_2$ and 
the fact that $S_2$ is flat over $\K$. So $S$ is flat over $S_1$ since 
$\otimes_\K S_2 =\otimes_{S_1} S_1 \otimes_{\K} S_2 =\otimes_{S_1} S$. 
 
First, show that $P_1 \otimes_\K S_2$ is a projective $S$-module if $P_1$ is a projective 
$S_1$-module. If $P_1$ is so, then $S_1^{\oplus n}=P_1 \oplus Q_1$, i.e., 
$P_1$ is a direct summand of a free $S_1$-module. 
Using the flatness of $\otimes_\K S_2$, it holds that 
$S^{\oplus n} =(P_1 \otimes_\K S_2) \oplus (Q_1 \otimes_\K S_2)$. Hence $P_1 \otimes_\K S_2$ is 
projective. Also, it is known that 
$$
\mbox{Hom}_{S_1} (M,S_1) \otimes_\K S_2 
\simeq \mbox{Hom}_{S} (S \cdot M,S)
$$
for any finitely generated $S_1$-module $M$ (see \cite{M} for example).
Hence 
$$
\mbox{Ext}_{S_1}^q (\Omega^p(\A_1,m_1),S_1) \otimes_\K S_2
\simeq \mbox{Ext}_S^q (S \cdot \Omega^p(\A_1,m_1),S).
$$ 
So it holds that 
$\pd_S S \cdot \Omega^p(\A_1,m_1) \ge \pd_{S_1} \Omega^p(\A_1,m_1)$. \owari
%Put $S_1=\K[x_1,\ldots,x_d]$. Then $x_{d+1},\ldots,x_\ell$ is a regular sequence of 
%$S\cdot \Omega^p(\A_1,m_1)$. Also, it is obvious that 
%$$
%S \cdot \Omega^p(\A_1,m_1)/(x_{d+1},\cdots,x_\ell) S \cdot \Omega^p(\A_1,m_1)=\Omega^p(\A_1,m_1).
%$$
%Thus $\depth_S S \cdot \Omega^p(\A_1,m_1)=\depth_{S_1} \Omega^p(\A_1,m_1)+\ell-d$. 
%Then Auslander-Buchsbaum shows that 
%\begin{eqnarray*}
%\pd_{S_1} \Omega^p(\A_1,m_1)&=&d-\depth_{S_1} \Omega^p(\A_1,m_1)\\
%&=&d-(\depth_S S \cdot \Omega^p(\A_1,m_1)-\ell+d)\\
%&=&\pd_S S \cdot \Omega^p(\A_1,m_1) \le \pd_S\Omega^p(\A,m),
%\end{eqnarray*}
%which completes the proof.\owari
\medskip

\noindent
\textbf{Proof of Theorem \ref{MCA}, continued}. 
Lemma \ref{pd} allows us to apply the induction hypothesis on dimensions to the essentialization of 
$\A_X$ and $(\A''_X,m_X)$ since they are also tame. Also, note that 
$\chi_{red}(\A''_X,m_X,t)$ is nothing but the characteristic polynomial 
of the essentialization of $(\A''_X,m_X)$. Hence 
$$
%|\chi_0({(\A)}_X^e,0)| 
b_i^X \ge |\chi_{red}(\A''_X,m_X,0)|=\sigma_i^X.
$$
Then local-to-global formula above shows that $b_i \ge \sigma_i\ (i=0,1,\ldots,\ell-2)$. 

Next show that 
$b_{\ell-1} \ge \sigma_{\ell-1}$. 
%By the arguments above combined with the 
%vanishing of $H^p(C^*)$ we saw before, it holds that 
%Note that we have not yet used the tameness assumption. From now on we use it.
By the assumption and Lemma \ref{lemma2}, it holds that $\pd_{S'} M^p \le p\ (p=0,1,\ldots,\ell-1)$. 
Hence Theorem \ref{TY72} combined with Proposition \ref{finite} shows 
that $H^p(M^*)=0\ (p \le \ell-2)$. Also, Theorem \ref{TY72} combined with 
Proposition \ref{finite} and the assumption that 
$(\A'',m)$ is tame show that 
$H^p(\Omega^*(\A'',m))=0\ (p \le \ell-2)$ for every 
generic $\eta_d$-complex. 
Hence the 
long exact sequence of cohomology of the sequence 
$$
0 \rightarrow M^p \rightarrow \Omega^p(\A'',m) \rightarrow C^p \rightarrow 0
$$
shows that 
$H^p(C^*)=0\ (0 \le p \le \ell-3)$ for every 
generic $\eta_d$-complex. 
By the arguments in \cite{Sc},
\begin{eqnarray*}
\chi_0(\A,t)&=& \sum_{p=0}^{\ell-1} \Poin (M^p,x) (t(1-x)-1)^p|_{x=1},\\
\chi(\A'',m,t)&=& \sum_{p=0}^{\ell-1} \Poin (\Omega^p(\A'',m),x) (t(1-x)-1)^p|_{x=1}.
\end{eqnarray*}
Hence 
$$
\chi_0(\A,t)-\chi(\A'',m,t)=-\sum_{p=0}^{\ell-1} \Poin (C^p,x)(t(1-x)-1)^p|_{x=1}.
$$
Now consider a generic $\eta_0$-complex combined with the cohomology vanishing 
in the above. Then 
\begin{eqnarray*}
\chi_0(\A,0)-\chi(\A'',m,0)&=&-\sum_{p=0}^{\ell-1} \Poin (C^p,x)(-1)^p|_{x=1}\\
&=&-\sum_{p=1}^{\ell-2} \Poin H^p (C^p)(-1)^p|_{x=1}\\
&=&(-1)^{\ell-1} \dim_\K H^{\ell-2}(C^*).
\end{eqnarray*}
Hence $b_{\ell-1}-\sigma_{\ell-1}=\dim_\K H^{\ell-2}(C^*) \ge 0$. 
To complete the proof, it suffices to show the following proposition. \owari
%%the previous sections show that 
%\begin{eqnarray*}
%\chi_0(\A,-1)-\chi(\A'',m,-1)&=&
%[\chi_0(\A,\displaystyle \frac{x-1}{x(1-x)})-\chi(\A'',m,\displaystyle \frac{x-1}{x(1-x)})]_{x=1}\\
%&=&-\sum_{p=0}^{\ell-1} \mbox{Poin}(C^p,x)(-x)^{-p}|_{x=1}\\
%&=&-\sum_{p=1}^{\ell-2} \mbox{Poin}(C^p,x)(-x)^{-p}|_{x=1}\\
%&=&-\sum_{p=1}^{\ell-2} \mbox{Poin}(H^p(C^*),x)(-x)^{-p}|_{x=1}\\
%&=&-\mbox{Poin}(H^{\ell-2}(C^*),x)(-x)^{\ell-2}|_{x=1}.
%\end{eqnarray*}
%Note that we use Proposition \ref{poly} to substitute several functions and so on. 
%Obviously a Poincar\'{e} 
%polynomial of a finite dimensional over $S$ has only non-negative coefficients, 
%which completes the proof. Finally, the non-negativity of $(-1)^{\ell-1} \chi(\A'',m,-1)$ 
%follows by the proposition below. \owari
\medskip

\begin{prop}
If $(\A'',m)$ is tame, then 
$\sigma_i \ge 0$. 
\label{nonnega1}
\end{prop}

\noindent
\textbf{Proof}. 
Use the similar argument to the proof of Theorem \ref{MCA}. Then it suffices to show that, 
by using localizations and Theorem \ref{ATW}, the Euler characteristic 
$\sum_{p=0}^{\ell-1} (-1)^p \mbox{Poin}(\Omega^p(\A'',m),x) |_{x=1}$ is not negative. 
Then the tameness condition completes the proof. \owari
%This value is nothing but the Euler characteristic of the $\eta_0$-complex 
%for a generic $\eta_0$ as in Proposition \ref{finite}. Then as in the proof of Theorem 
%\ref{MCA}, the tameness condition shows that the value is nothing but 
%the dimension of the highest cohomology group. Here, also we used the fact that 
%$(\A_X,m_X)$ is tame if $(\A,m)$ is tame. 
%\owari

\begin{cor}
In the notation of Theorem \ref{MCA}, 
$(-1)^{\ell-1} \chi(\A'',m,-1) \ge 0$. 
%Hence 
%Theorem \ref{MCA} gives a non-negative lower bound of chambers of 
%the deconing of $\A$.
\label{p}
\end{cor}

\noindent
\textbf{Proof}. 
Proposition \ref{p} and 
$$
(-1)^{\ell-1} \chi_0(\A'',m,-1)=
\sum_{i=0}^{\ell-1} \sigma_i
$$
completes the proof. \owari
%By using the same arguments as in the proof of Theorem \ref{MCA}, 
%it holds that 
%\begin{eqnarray*}
%\chi(\A'',m,-1)&=&
%\sum_{p=0}^{\ell-1} \mbox{Poin}(H^p(\Omega(\A'',m)),x)(-x)^{-p}|_{x=1}\\
%&=&\mbox{dim}_\K H^{\ell-1}(\Omega^*	(\A'',m)))(-1)^{\ell-1}, 
%\end{eqnarray*}
%which completes the proof. \owari
\medskip

\noindent
\textbf{Proof of Corollary \ref{MCA2}}. 
The first statement follows immediately from Theorem \ref{MCA} and 
Corollary \ref{p}. Let us prove the second statement. By \cite{Z} and the 
factorization in \cite{ATW}, it is easy to see that the freeness implies 
MCA. Assume that $\A$ is MCA. Then Theorem \ref{MCA} implies that 
$\chi_0(\A,t)=\chi(\A'',m,t)$. Then Theorem \ref{AYfree} 
%or the 
%freeness criterion in \cite{AY} 
completes the proof. \owari

%Also, note that the lower bound in Theorem \ref{MCA} is non-negative.

\begin{cor}
Assume that $\A$ is a $4$-arrangement. Then the statement in Theorem \ref{MCA} 
holds true if $\A$ is tame.
\end{cor}

\noindent
\textbf{Proof}. 
Since $\Omega^p(\A'',m)$ is reflexive, Auslander-Buchsbaum formula combined with 
$\depth_{S'} \Omega^p(\A'',m) \ge 2$ completes the proof. \owari
%Note that every $3$-multiarrangement is locally free. Then Corollary 5.4 in 
%\cite{MS} confirms that $(\A'',m)$ is tame. \owari
\medskip

%By these results we have interesting corollaries.
%
%\begin{cor}
%Assume the whole tameness and $(\A'',m)$ is free. 
%Then $\A$ is a minimal chamber arrangement if and only if 
%$\A$ is a minimal bounded chamber arrangement.
%\end{cor}

Apparently $2$ and $3$-(multi)arrangements are tame. So as a corollary of 
Theorem \ref{MCA} we can prove Yoshinaga's criterion. 

\begin{cor}[\cite{Y}, Theorem 3.2]
A $3$-arrangement $\A$ is free if and only if 
it is a minimal chamber arrangement.
% in the above sense.
\label{Y2}
\end{cor}

\noindent
\textbf{Proof}. 
Since $2$ and $3$-multiarrangements are tame, we can use Theorem \ref{MCA} and 
Corollary \ref{MCA2}. 
Assume that $\A$ is a minimal chamber arrangement. Since 
$C^0=C^2=0$, the complex is $ 0 \rightarrow C^1 \rightarrow 0$. The minimality of 
chambers implies that $H^1(C^*)=0$, which is nothing but $C^1=0$. Then 
the result in \cite{Z} implies that $\A$ is free. 
\owari
%the only 
%suviving cohomology is the middle one, and theorems say that the middle cohomology 
%also vanishes, we have the corollary. \owari
\medskip

Theorem \ref{MCA} and Corollary \ref{MCA2} give us a direction of 
the research on the relation between free arrangements and MCA as follows:

\begin{problem}
Do Theorem \ref{MCA} and Corollary \ref{MCA2} hold 
true without the assumption of tameness?
\label{pb}
\end{problem}

Problem \ref{pb} is a very natural one. When 
we construct a free arrangement by the addition theorem, we usually, or 
empirically,  
add hyperplanes in such a way that the new arrangements have the smallest 
chambers among all the other choices (though the addition of this type  
does not always work well!). This choice of the additions is jutified when $\ell =3$ by \cite{Y}. 
If Problem \ref{pb} is 
true, then we can obtain a better generalization of Yoshinaga's 
criterion. If that is not true, then the tameness condition becomes more 
important, and essential condition which connectes algebra and 
geometry of hyperplane arrangements.

\medskip

%\section{Related results}
%
%\begin{theorem}
%If $\A$ and $(\A'',m)$ are both tame, then 
%$\pi_0(\A,t)-\pi(\A'',m,t) \ge 0$.
%\label{nonnega2}
%\end{theorem}
%
%\noindent
%\textbf{Proof}. 
%Again use the local-to-global formula. Then it holds that, to compare 
%each coefficients, it suffices to compare that of the top term of 
%$\pi(\A''_Y,m_Y)\ (Y \in L(\A'')$ and that of 
%$\pi_0(\B_Y,t)$, where $\B_Y \subset \A$ is the extension of $(\A_Y,m_Y)$ in $\A$. 
%For this correspondence, see \cite{AY} for details. Then again 
%the (local) tameness and the Euler characteristic argument completes the proof. \owari
%\medskip
%
%\begin{cor}
%Assume that $\A$ is tame and $(\A'',m)$ is free. Then 
%$\A$ is free if and only if $\A$ is MCA.
%\label{nonnega3}
%\end{cor}
%
%\noindent
%\textbf{Proof}. The ``only if'' part is trivial. 
%Assume that $\A$ is MCA. Then Theorem \ref{nonnega2} shows that 
%$\pi_0(\A,t)=\pi(\A'',m,t)$. Then Theorem \ref{Schulze} completes the proof. \owari

 \vspace{5mm}

\noindent
Takuro Abe\\
Department of Mechanical Engineering and Science\\
Kyoto University\\
Yoshida Honmachi, Sakyo-Ku, Kyoto 6068501, Japan\\
abe.takuro.4c@kyoto-u.ac.jp

\end{document}